\input amstex
\documentstyle{amsppt}
\magnification=\magstep1                        
\hsize6.5truein\vsize8.9truein                  
\NoRunningHeads
\loadeusm

\magnification=\magstep1                        
\hsize6.5truein\vsize8.9truein                  
\NoRunningHeads
\loadeusm

\document
\topmatter

\title
The asymptotic value of the Mahler measure of the Rudin-Shapiro polynomials
\endtitle

\rightheadtext{the Mahler measure of the Rudin-Shapiro polynomials}

\author Tam\'as Erd\'elyi
\endauthor

\address Department of Mathematics, Texas A\&M University,
College Station, Texas 77843, College Station, Texas 77843 (T. Erd\'elyi) \endaddress

\thanks {{\it 2010 Mathematics Subject Classifications.} 11C08, 41A17, 26C10, 30C15}
\endthanks

\keywords
polynomial inequalities, Mahler measure, Rudin-Shapiro polynomials, zeros of polynomials
\endkeywords

\date September 9, 2017 
\enddate

\email terdelyi\@math.tamu.edu
\endemail

\abstract
In signal processing the Rudin-Shapiro polynomials have good autocorrelation properties
and their values on the unit circle are small. Binary sequences with low autocorrelation
coefficients are of interest in radar, sonar, and communication systems.
In this paper we show that the Mahler measure of the Rudin-Shapiro polynomials of degree 
$n-1$ with $n=2^k$ is asymptotically $(2n/e)^{1/2}$, as it was conjectured by B. Saffari in 
1985. Our approach is based heavily on the Saffari and Montgomery conjectures proved 
recently by B. Rodgers.  
\endabstract

\endtopmatter

\head 1. Introduction and Notation \endhead

Let $D := \{z \in {\Bbb C}: |z| < 1\}$ denote the open unit disk of the complex plane.
Let $\partial D :=  \{z \in {\Bbb C}: |z| = 1\}$ denote the unit circle of the complex plane.
The Mahler measure $M_{0}(f)$ is defined for bounded measurable functions $f$ on $\partial D$ by 
$$M_{0}(f) := \exp\left(\frac{1}{2\pi} \int_{0}^{2\pi}{\log|f(e^{it})|\,dt} \right)\,.$$
It is well known, see [HL-52], for instance, that
$$M_{0}(f) = \lim_{q \rightarrow 0+}{M_{q}(f)}\,,$$
where
$$M_{q}(f) := \left( \frac{1}{2\pi} \int_{0}^{2\pi}{\left| f(e^{it}) \right|^q\,dt} \right)^{1/q}\,, 
\qquad q > 0\,.$$
It is also well known that for a function $f$ continuous on $\partial D$ we have 
$$M_{\infty}(f) := \max_{t \in [0,2\pi]}{|f(e^{it})|} = \lim_{q \rightarrow \infty}{M_{q}(f)}\,.$$
It is a simple consequence of the Jensen formula that
$$M_0(f) = |c| \prod_{j=1}^n{\max\{1,|z_j|\}}$$
for every polynomial of the form
$$f(z) = c\prod_{j=1}^n{(z-z_j)}\,, \qquad c,z_j \in {\Bbb C}\,.$$
See [BE-95, p. 271] or [B-02, p. 3], for instance.

Let ${\Cal P}_n^c$ be the set of all algebraic polynomials of degree at most $n$ with complex coefficients. 
Let ${\Cal T}_n$ be the set of all real (that is, real-valued on the real line) trigonometric polynomials 
of degree at most $n$ with real coefficients. Finding polynomials with suitably restricted coefficients and 
maximal Mahler measure has interested many authors. The classes
$${\Cal L}_n := \left\{ f: \enskip f(z) = \sum_{j=0}^{n}{a_jz^j}\,, \quad a_j \in \{-1,1\} \right\}$$
of Littlewood polynomials and the classes
$${\Cal K}_n := \left\{ f: \enskip f(z) = \sum_{j=0}^{n}{a_jz^j}\,, \quad a_j \in {\Bbb C}, \enskip |a_j| =1 \right\}$$
of unimodular polynomials are two of the most important classes considered.
Observe that ${\Cal L}_n \subset {\Cal K}_n$ and
$$M_0(f) \leq M_2(f) = \sqrt{n+1}$$
for every $f \in {\Cal K}_n$.
Beller and Newman [BN-73] constructed unimodular polynomials $f_n \in {\Cal K}_n$  whose
Mahler measure $M_0(f_n)$ is at least $\sqrt{n}-c/\log n$.

Section 4 of [B-02] is devoted to the study of Rudin-Shapiro polynomials. 
Littlewood asked if there were polynomials $f_{n_k} \in {\Cal L}_{n_k}$ satisfying  
$$c_1 \sqrt{n_k+1}  \leq |f_{n_k}(z)| \leq c_2 \sqrt{n_k+1}\,, \qquad z \in \partial D\,,$$
with some absolute constants $c_1 > 0$ and $c_2 > 0$, see [B-02, p. 27] for a reference 
to this problem of Littlewood.
To satisfy just the lower bound, by itself, seems very hard, and no such sequence $(f_{n_k})$  
of Littlewood polynomials $f_{n_k} \in {\Cal L}_{n_k}$ is known. A sequence of Littlewood polynomials 
that satisfies just the upper bound is given by the Rudin-Shapiro polynomials. The Rudin-Shapiro 
polynomials appear in Harold Shapiro's 1951 thesis [S-51] at MIT and are sometimes called just 
Shapiro polynomials. They also arise independently in Golay's paper [G-51]. They are 
remarkably simple to construct and are a rich source of counterexamples to possible 
conjectures.

The Rudin-Shapiro polynomials are defined recursively as follows:
$$\split P_0(z) & :=1\,, \qquad Q_0(z) := 1\,, \cr 
P_{k+1}(z) & := P_k(z) + z^{2^k}Q_k(z)\,, \cr
Q_{k+1}(z) & := P_k(z) - z^{2^k}Q_k(z)\,, \cr \endsplit$$
for $k=0,1,2,\ldots\,.$ Note that both $P_k$ and $Q_k$ are polynomials of degree $n-1$ with $n := 2^k$
having each of their coefficients in $\{-1,1\}$.
In signal processing, the Rudin-Shapiro polynomials have good autocorrelation properties 
and their values on the unit circle are small. Binary sequences with low autocorrelation 
coefficients are of interest in radar, sonar, and communication systems.

It is well known and easy to check by using the parallelogram law that
$$|P_{k+1}(z)|^2 + |Q_{k+1}(z)|^2 = 2(|P_k(z)|^2 + |Q_k(z)|^2)\,, \qquad z \in \partial D\,.$$
Hence
$$|P_k(z)|^2 + |Q_k(z)|^2 = 2^{k+1} = 2n\,, \qquad z \in \partial D\,. \tag 1.1$$
It is also well known (see Section 4 of [B-02], for instance), that
$$|Q_k(z)| = |P_k(-z)|\,, \qquad z \in \partial D\,. \tag 1.2$$
P. Borwein's book [B-02] presents a few more basic results on the Rudin-Shapiro
polynomials. Various properties of the Rudin-Shapiro polynomials are discussed in
[B-73] and [BL-76]. Obviously $M_2(P_k) = 2^{k/2}$ by the Parseval formula. In 1968 Littlewood
[L-68] evaluated $M_4(P_k)$ and found that $M_4(P_k) \sim (4^{n+1}/3)^{1/4}$.
The $M_4$ norm of Rudin-Shapiro like polynomials on $\partial D$ are studied in [BM-00]. 
P. Borwein and Lockhart [BL-01] investigated the asymptotic behavior of the mean
value of normalized $M_q$ norms of Littlewood polynomials for arbitrary $q > 0$.
They proved that
$$\lim_{n \rightarrow \infty} {\frac{1}{2^{n+1}} \, \sum_{f \in {\Cal L}_n}
{\frac{(M_q(f))^q}{n^{q/2}}}}= \Gamma \left( 1+ \frac q2 \right)\,.$$
In [C-15c] we proved that
$$\lim_{n \rightarrow \infty} {\frac{1}{2^{n+1}} \sum_{f \in {\Cal L}_n}
{\frac{M_q(f)}{n^{1/2}}}} = \left( \Gamma \left( 1+ \frac q2 \right) \right)^{1/q}$$
for every $q > 0$. In [CE-15c] we showed also that
$$\lim_{n \rightarrow \infty} {\frac{1}{2^{n+1}} \sum_{f \in {\Cal L}_n}
{\frac{M_0(f)}{n^{1/2}}}} = e^{-\gamma/2}\,,$$
where
$$\gamma := \lim_{n \rightarrow \infty}{\left( \sum_{k=1}^n{\frac 1k - \log n} \right)} = 0.577215 \ldots$$
is the Euler constant and $e^{-\gamma/2} = 0.749306\ldots$.
These are analogues of the results proved earlier by Choi and Mossinghoff
[CM-11] for polynomials in ${\Cal K}_n$.

In 1980 Saffari conjectured the following.

\proclaim{Conjecture 1.1}
Let $P_k$ and $Q_k$ be the Rudin-Shapiro polynomials of degree $n-1$ with $n := 2^k$. 
We have 
$$M_q(P_k) = M_q(Q_k) \sim \frac{2^{(k+1)/2}}{(q/2+1)^{1/q}}$$
for all real exponents $q > 0$. Equivalently, we have 
$$\split & \lim_{n \rightarrow \infty} 
m{\left(\left\{t \in K: \left| \frac{P_k(e^{it})}{\sqrt{2^{k+1}}} \right|^2 \in [\alpha,\beta] \right\}\right)} \cr 
= \, & \lim_{n \rightarrow \infty}
m{\left(\left\{t \in K: \left| \frac{Q_k(e^{it})}{\sqrt{2^{k+1}}} \right|^2 \in [\alpha,\beta] \right\}\right)} 
= \beta - \alpha \cr \endsplit$$
whenever $0 \leq \alpha < \beta \leq 1$. 
\endproclaim

This conjecture was proved for all even values of $q \leq 52$ by Doche [D-05]
and Doche and Habsieger [DH-04]. Recently B. Rodgers [R-16] proved Saffari's Conjecture 1.1 
for all $q > 0$. See also [EZ-17]. An extension of Saffari's conjecture is Montgomery's conjecture below.  

\proclaim{Conjecture 1.2}
Let $P_k$ and $Q_k$ be the Rudin-Shapiro polynomials of degree $n-1$ with $n := 2^k$.
We have
$$\split & \lim_{n \rightarrow \infty} 
m{\left(\left\{t \in K: \frac{P_k(e^{it})}{\sqrt{2^{k+1}}} \in E \right\}\right)} \cr
= \, & \lim_{n \rightarrow \infty}
m{\left(\left\{t \in K: \frac{Q_k(e^{it})}{\sqrt{2^{k+1}}} \in E \right\}\right)} 
= \frac{1}{\pi} \, m(E) \cr \endsplit$$
for any measurable set $E \subset D := \{z \in {\Bbb C}: |z| < 1\}\,.$
\endproclaim

B. Rodgers [R-16] proved Montgomery's Conjecture 1.2 as well. 

Despite the simplicity of their definitions not much is known
about the Rudin-Shapiro polynomials. It has been shown in [E-16] fairly recently that 
the Mahler measure ($M_0$ norm) and the $M_\infty$ norm of the Rudin-Shapiro polynomials 
$P_k$ and $Q_k$ of degree $n-1$ with $n := 2^k$ on the unit circle of the complex plane 
have the same size, that is, the Mahler measure of the Rudin-Shapiro polynomials   
of degree $n-1$ with $n := 2^k$ is bounded from below by $cn^{1/2}$, where $c > 0$ is 
an absolute constant.

It is shown in this paper that the Mahler measure of the Rudin-Shapiro polynomials $P_k$ and $Q_k$ 
of degree $n-1=2^k-1$ is asymptotically $(2n/e)^{1/2}$, as it was conjectured by B. Saffari 
in 1985. Note that $(2/e)^{1/2} = 0.85776388496\ldots$ is larger than $e^{-\gamma/2} = 0.749306\ldots$ 
in the average Mahler measure result for the class of Littlewood polynomials ${\Cal L}_n$ we mentioned 
before.  

\head 2. New Result \endhead

Let $P_k$ and $Q_k$ be the Rudin-Shapiro polynomials of degree $n-1$ with $n := 2^k$.

\proclaim{Theorem 2.1}
We have 
$$\lim_{n \rightarrow \infty}{\frac{M_0(P_k)}{n^{1/2}}} = \lim_{n \rightarrow \infty}{\frac{M_0(Q_k)}{n^{1/2}}} 
= \left( \frac 2e \right)^{1/2}\,.$$  
\endproclaim

\head 3. Lemmas \endhead 

Let $D(a,r) := \{z \in {\Bbb C}: |z-a| < r\}$ denote the open disk of the complex plane  centered at 
$a \in {\Bbb C}$ of radius $r > 0$.

To prove our theorem we need some lemmas. Our first lemma states Jensen's formula. Its proof may be 
found in most of the complex analysis textbooks.  

\proclaim{Lemma 3.1} Suppose
$h$ is a nonnegative integer and
$$f(z) = \sum^\infty_{k=h} {c_k (z-z_0)^k}\,, \qquad c_h \neq 0 $$
is analytic on the closure of the disk $D(z_0,r)$. Let $a_1, a_2, \ldots, a_m$ denote  
the zeros of $f$ in $D(z_0,r) \setminus \{z_0\}$, where
each zero is listed as many times as its multiplicity. We have
$$\log|c_h| + h \log r + \sum_{k=1}^m \log {r \over{|a_k-z_0|}}= {1 \over{2 \pi}}
\int_0^{2 \pi} {\log |f(z_0 + re^{i \theta})| \, d\theta}\,.$$
\endproclaim

\proclaim{Lemma 3.2}
There exists a constant $c_1$ depending only on $c_2 > 0$ such that every polynomial $P \in {\Cal P}_n^c$    
has at most $c_1(nr+1)$ zeros in any open disk $D(z_0,r)$
with $z_0 \in \partial D$ and 
$$|P(z_0)| \geq c_2M_{\infty}(P) \tag 3.1$$
\endproclaim

\demo{Proof of Lemma 3.2}
Without loss of generality we may assume that $z_0 := 1$ and 
$$n^{-1} \leq r \leq 1\,.$$ 
Indeed, the case $0 < r < n^{-1}$ follows from the case $r = n^{-1}$, and the case $r > 1$ is obvious.
Let $P$ be a polynomial of the form given in the lemma. 
A well known polynomial inequality observed by Bernstein states that  
$$|P(\zeta)| \leq \max\{1,|\zeta|^n\} \max_{z \in \partial D}{|P(z)|} \tag 3.2$$
for any polynomials $P \in {\Cal P}_n^c$ and for any $\zeta \in {\Bbb C}$.
This is a simple consequence of the Maximum Principle, see [BE-95, p. 239], for instance.
Using (3.2) we can deduce that  
$$\log |P(z)| \leq \log((1+2r)^n M_\infty(P)) \leq \log M_\infty(P) + 2nr\,, \qquad |z| \leq 1 + 2r\,. \tag 3.3$$
Let $m$ denote the number of zeros of $P$ in the open disk $D(z_0,r)$.
Using Lemma 3.1 with the disk $D(z_0,2r)$ and $h=0$, then using (3.1) and (3.3), we obtain
$$\log c_2 + \log M_\infty(P) + m \log 2 \leq \log |P(z_0)| + m \log 2 \leq \frac{1}{2\pi}\,2\pi (\log M_\infty + 2nr)\,.$$
This, together with $n^{-1} \leq r \leq 1$, implies 
$\log c_2 + m \log 2 \leq 2nr$, and the lemma follows.
\qed \enddemo

Our next lemma is stated as Lemma 3.5 in [E-16], where its proof may also be found.

\proclaim{Lemma 3.3}
If $P_k$ and $Q_k$ are the $k$-th Rudin-Shapiro polynomials of degree $n-1$ with $n:=2^k$, 
$\gamma := \sin^2(\pi/8)$, and
$$z_j := e^{it_j}\,, \quad t_j := \frac{2\pi j}{n}\,, \qquad j \in {\Bbb Z}\,,$$
then
$$\max \{|P_k(z_j)|^2,|P_k(z_{j+r})|^2\} \geq \gamma 2^{k+1} = 2\gamma n\,, \quad r \in \{-1,1\}\,,$$
for every $j=2u$, $u \in {\Bbb Z}$.
\endproclaim

By Lemma 3.3, for every $n = 2^k$  there are
$$0 \leq \tau_1 < \tau_2 < \cdots < \tau_m < \tau_{m+1} := \tau_1 + 2\pi$$
such that 
$$\tau_j - \tau_{j-1} = \frac{2\pi l}{n}\,, \qquad l \in \{1,2\}\,,$$
and with 
$$a_j := e^{i\tau_j}, \qquad j=1,2,\ldots,m+1, \tag 3.4$$
we have
$$|P_k(a_j)|^2 \geq 2\gamma n\,, \qquad j=1,2,\ldots,m+1\,. \tag 3.5$$
(Moreover, each $a_j$ is an $n$-th root of unity.) For the sake of brevity let $R_n \in {\Cal T}_n$ 
be defined by 
$$R_n(t) := |P_k(e^{it})|^2, \qquad n=2^k\,.$$

Using the above notation we formulate the following observation.

\proclaim{Lemma 3.4}
There is an absolute constant $c_3 > 0$ such that
$$\mu := \left|\left \{j \in \{2,3,\ldots,m+1\}: \min_{t \in [\tau_{j-1},\tau_j]}{R_n(t)} \leq \varepsilon \right \}\right| 
\leq c_3n\varepsilon^{1/2}$$
for every sufficiently large $n=2^k \geq n_\varepsilon$, $k=1,2,\ldots$, and $\varepsilon > 0$.
\endproclaim

To prove Lemma 3.4 we need a consequence of the so-called Bernstein-Szeg\H o inequality 
formulated by our next lemma. For its proof see [BE-95, p. 232], for instance. 

\proclaim{Lemma 3.5}
We have 
$$S^\prime(t)^2 + n^2S(t)^2 \leq n^2 \max_{\tau \in {\Bbb R}} S(\tau)^2$$
for every $S \in {\Cal T}_n$.
\endproclaim

\proclaim{Lemma 3.6}
We have
$$|R_n^\prime(t)| \leq n^{3/2} \sqrt{2R_n(t)}\,, \qquad t \in {\Bbb R}\,.$$
\endproclaim

\demo{Proof of Lemma 3.6}
Let $S \in {\Cal T}_n$ be defined by $S(t) := |P_k(e^{it})|^2-n = R_n(t)-n$.
Observe that (1.1) implies that
$$\max_{\tau \in {\Bbb R}}{|S(\tau)|} \leq n\,.$$
Combining this with Lemma 3.5 implies that
$$\split |R_n^\prime(t)| & = |S^\prime(t)| = n \sqrt{n^2-S(t)^2} \leq n \sqrt{n^2-(R_n(t)-n)^2} + n \sqrt{2nR_n(t) - R_n(t)^2)} \cr 
& \leq n\sqrt{2nR_n(t)}\,. \cr \endsplit$$         
\qed \enddemo

Now we are ready to prove Lemma 3.4.

\demo{Proof of Lemma 3.4}
Let $j \in \{2,3,\ldots,m+1\}$ be such that 
$$\min_{t \in [\tau_{j-1},\tau_j]}{R_n(t)} \leq \varepsilon \,.$$
Using the notation of Lemma 3.3, without loss of generality we may assume that $0 < \varepsilon < \gamma$. 
iBy recalling (3.4) and (3.5) there are $\tau_{j-1} \leq \alpha_j < \beta_j \leq \tau_j$ such that
$$R_n(\alpha_j) = \varepsilon n\,, \qquad R_n(\beta_j) = 2\varepsilon n\,,$$
and
$$R_n(t) \leq 2\varepsilon n\,, \qquad t \in [\alpha_j,\beta_j]\,.$$ 
Then, by the Mean Value Theorem there is $\xi_j \in (\alpha_j,\beta_j)$ such that 
$$\varepsilon n = R_n(\beta_j) - R_n(\alpha_j) = (\beta_j - \alpha_j)R_n^\prime(\xi_j)\,,$$
and hence by Lemma 3.5 we obtain
$$\split \varepsilon n & = (\beta_j - \alpha_j)R_n^\prime(\xi_j) \leq  (\beta_j - \alpha_j)n^{3/2} \sqrt{2R_n(\xi_j)} \cr
& \leq (\beta_j - \alpha_j)n^{3/2} \sqrt{4\varepsilon n}\,, \cr \endsplit$$
that is
$$\beta_j - \alpha_j \geq \frac{\varepsilon^{1/2}}{2n}\,.$$  
Hence, on one hand, 
$$m{\left(\left\{t \in K: \frac{R_n(t)}{n} \in [0,2\varepsilon] \right\}\right)} = 
m{\left(\left\{t \in K: \left| \frac{P_k(e^{it})}{\sqrt{2^{k+1}}} \right|^2 \in [0,\varepsilon] \right\}\right)}
\geq \frac{\mu \varepsilon^{1/2}}{2n}\,.$$
On the other hand, by Conjecture 1.1 proved by B. Rodgers there is an absolute constant $c_3/2 > 0$ such that 
$$m{\left(\left\{t \in K: \frac{R_n(t)}{n} \in [0,2\varepsilon] \right\}\right)} = 
m{\left(\left\{t \in K: \left| \frac{P_k(e^{it})}{\sqrt{2^{k+1}}} \right|^2 \in [0,\varepsilon] \right\}\right)}
\leq (c_3/2) \varepsilon$$
for every sufficiently large $n \geq n_\varepsilon$.
Combining the last two inequalities we obtain
$$\mu \leq c_3n\varepsilon^{1/2}$$
for every sufficiently large $n \geq n_\varepsilon$.
\qed \enddemo

We introduce the notation 
$$A_{n,\varepsilon} := \left \{t \in K: \frac{R_n(t)}{2n} \geq \varepsilon \right \}$$
and
$$B_{n,\varepsilon} := K \setminus A_{n,\varepsilon} = \left \{t \in K: \frac{R_n(t)}{2n} < \varepsilon \right \}\,.$$

Our next lemma is an immediate consequence of Conjecture 1.1 proved by B. Rodgers.

\proclaim{Lemma 3.7}
Let $\varepsilon \in (0,1)$ be fixed. We have
$$\lim_{n \rightarrow \infty}{\frac{1}{2\pi} \int_{A_{n,\varepsilon}}{\log \frac{R_n(t)}{2n}\, dt}} = \int_\varepsilon^1{\log x\, dx}\,.$$ 
\endproclaim

\demo{Proof of Lemma 3.7}
Let 
$$F_\varepsilon(x) := \cases
\log x, \quad &\text{if \enskip} x \in [\varepsilon,1]\,,
\\
0, \quad &\text{if \enskip} x \in [0,\varepsilon)\,.
\endcases$$ 
By using the Weierstrass Approximation Theorem, it is easy to see that $F_\varepsilon$ can be approximated by polynomials in $L_1[0,1]$ norm, and hence 
the lemma follows from Conjecture 1.1 proved by B. Rodgers in a standard fashion. We omit the details of 
this routine argument. 
\qed \enddemo

The above lemma will be coupled with the following inequality.

\proclaim{Lemma 3.8}
Let $\varepsilon \in (0,1)$ be fixed. There is an absolute constant $c_4 > 0$ such that
$$\frac{1}{2\pi} \int_{B_{n,\varepsilon}}{\log \frac{R_n(t)}{2n}\, dt} \geq -c_4\varepsilon^{1/2}$$
for every sufficiently large $n \geq n_\varepsilon$.
\endproclaim

To prove Lemma 3.8 we need a few other lemmas.  

\proclaim{Lemma 3.9}
Let $f$ be a twice differentiable function on $[a,b]$. There is a $\xi \in [a,b]$ such that 
$$\int_a^b{f(t)\,dt} - \frac 12 (f(a) + f(b))(b-a) = -\frac{(b-a)^3}{12}{f^{\prime\prime}(\xi)}\,.$$
\endproclaim

This is the formula for the error term in the trapezoid rule. Its proof may be found in 
various calculus textbooks discussing numerical integration.

Let $w_j \in  {\Bbb C}$, $j=1,2,\ldots,n-1$, denote the zeros of $P_k$. So we have
$$\log{\frac{R_n(t)}{2n}} = \log{\frac{|P_k(e^{it})|^2}{n}} = \sum_{j=1}^{n-1}{\log{|e^{it}-w_j|}} - \log n\,.$$
It is a simple well known fact that $P_k \in {\Cal L}_{n-1}$ implies that 
$$1/2 \leq |w_j| \leq 2, \qquad j=1,2,\ldots,n-1\,.$$
Associated with $w \in {\Bbb C}$ we introduce $\phi \in [0,2\pi)$ uniquely defined by $w = |w|e^{i\phi}$. 
For the sake of brevity let
$$g_w(t) :=  \log{|e^{it}-w|} = \log{|e^{it}-|w|e^{i\phi}|}\,.$$
Simple calculations show that
$$g_w(t) = \frac 12\log(1 + |w|^2 - 2|w|\cos(t-\phi))\,,$$
$$g_w^\prime(t) = \frac{|w|\sin(t-\phi)}{|e^{it}-w|^2}\,,$$
and
$$g_w^{\prime\prime}(t) = \frac{|w|\cos(t-\phi)}{|e^{it}-w|^4} - \frac{2|w|^2\sin^2(t-\phi)}{|e^{it}-w|^4}\,.$$
The inequality of the following lemma is immediate.

\proclaim{Lemma 3.10}
There is an absolute constant $c_5 > 0$ such that 
$$|g_w^{\prime\prime}(t)| \leq \frac{c_5}{|e^{it}-w|^2}$$
for every $t \in {\Bbb R}$ and $w \in {\Bbb C}$ with $|w| \leq 2$.
\endproclaim

Combining Lemmas 3.9 and 3.10 we get the following.

\proclaim{Lemma 3.11}
Let $a_j = e^{i\tau_j}, j=1,2,\ldots m+1,$ be as before (defined after Lemma 3.3). 
There are $\xi_j \in [\tau_{j-1},\tau_j]$ and an absolute constant $c_6 > 0$ such that 
$$\int_{\tau_{j-1}}^{\tau_j}{g_w(t) \, dt} - \frac 12 (g_w(\tau_j) - g_w(\tau_{j-1})(\tau_j-\tau_{j-1}) \geq 
-\frac{c_6}{n^3|e^{i\xi_j}-w|^2}$$
for every $t \in {\Bbb R}$ and $w \in {\Bbb C}$ with $|w| \leq 2$.
\endproclaim

We will also need an estimate better than the one given in Lemma 3.11 in the case when 
$w \in {\Bbb C}$ is close to $a_j := e^{i\tau_j}$.

\proclaim{Lemma 3.12}
Let $a_j = e^{i\tau_j}, j=2,3,\ldots m+1$ be as before (defined after Lemma 3.3). 
There is an absolute constant $c_7 > 0$ such that 
$$\int_{\tau_{j-1}}^{\tau_j}{g_w(t) \, dt} - \frac 12 (g_w(\tau_j) - g_w(\tau_{j-1})(\tau_j-\tau_{j-1}) \geq -\frac{c_7}{n}$$
for every $t \in {\Bbb R}$ and $w \in {\Bbb C}$ such that $|w-a_j| \leq 8\pi/n$.
\endproclaim

To prove Lemma 3.12 we need the following observation.

\proclaim{Lemma 3.13}
Let $a_j = e^{i\tau_j}, j=1,2,\ldots m+1,$ be as before (defined after Lemma 3.3).
There is an absolute constant $c_8 > 0$ such that $|a_j-w| \geq c_8/n$ for every $w \in {\Bbb C}$ 
for which $P(w)=0$.
\endproclaim

\demo{Proof of Lemma 3.13}
The proof is a routine combination of (1.1), Lemma 3.3, and a couple of Bernstein's inequalities.
One of Bernstein's polynomial inequalities asserts that 
$$\max_{z \in \partial D}{|P^{\prime}(z)|} \leq n \, \max_{z \in \partial D}{|P(z)|} \tag 3.6$$
for any polynomials $P \in {\Cal P}_n^c$. See [BE-95, p. 232], for instance. Another polynomial inequality of 
Bernstein we need in this proof is (3.2). 

Suppose $w \in {\Bbb C}$, $|w-a_j| \leq c/n$, and $P(w)=0$. Let $\Gamma$ be the line segment connecting 
$a_j = e^{i\tau_j}$ and $w$. We have
$$\split (2\gamma)^{1/2}n^{1/2} \leq & |P_k(a_j)| = |P_k(a_j) - P_k(w)| 
= \left| \int_{\Gamma} {P_k^\prime(z) \, dz} \right| \cr 
\leq & \int_{\Gamma}{|P_k^\prime(z)| \, |dz|}\,. \cr \endsplit$$ 
Hence there is a $\zeta \in \Gamma$ such that
$$|P_k^\prime(\zeta)| \cdot |a_j-w| \geq (2\gamma)^{1/2}n^{1/2}\,.$$ 
Combining this  with (3.6) and $|\zeta-a_j| \leq |w-a_j| \leq c/n$, we obtain
$$|P_k^\prime(\zeta)| \geq \frac{(2\gamma)^{1/2}n^{1/2}}{c/n} = \frac{(2\gamma)^{1/2}}{c}n^{3/2}\,. \tag 3.7$$
On the other hand, combining (3.2), (3.6), and (1.1), we obtain
$$\split |P_k^\prime(\zeta)| \leq & \left( \max\{1,|\zeta|^{n-1}\} \right) \left( \max_{z \in \partial D}{|P_k^\prime(z)|} \right) \cr  
\leq & \left( \max\{1,|\zeta|^{n-1}\} \right) \left(n\max_{z \in \partial D}{|P_k(z)|}\right) \leq \left(1 + \frac cn \right)^n n (2n)^{1/2} \cr
\leq & e^c \sqrt{2} n^{3/2} \cr \endsplit$$ 
Combining this with (3.7), we get $\gamma^{1/2} \leq ce^c$. 
\qed \enddemo 

\demo{Proof of Lemma 3.12}
Observe that Lemma 3.13 implies that there is an absolute constant $c_8 > 0$ such that
$$\frac 12 (\log|e^{i\tau_j}-w| + \log|e^{i\tau_{j-1}}-w|) \leq \log(c_8/n) = \log c_8 - \log n\,. \tag 3.8$$
Now we show that there is an absolute constant $c_9 > 0$ such that 
$$\int_{\tau_{j-1}}^{\tau_j}{\log|e^{it}-w| \, dt} \geq (\tau_j-\tau_{j-1})(c_9-\log n )\,. \tag 3.9$$
To see this let $w = |w|e^{i\phi}$. We have
$$|e^{it}-w| = |e^{it}-|w|e^{i\phi}| \geq |e^{it}-e^{i\phi}| = 2\sin \left|\frac{t-\phi}{2}\right|  
\geq 2 \, \frac{2}{\pi} \frac{|t-\phi|}{2} = \frac{2}{\pi} \, |t-\phi|$$
whenever $|t-\phi| \leq \pi$.
Hence, if $$\phi \in \left[\tau_{j-1} + \frac{c_8}{2n},\tau_j - \frac{c_8}{2n}\right]\,,$$ 
then 
$$\split \int_{\tau_{j-1}}^{\tau_j}{\log|e^{it}-w| \, dt} & \geq 
\int_{\tau_{j-1}}^{\tau_j}{\log \left( \frac{2}{\pi} \, |t-\phi| \right) \, dt} \cr
& = \int_{\tau_{j-1}}^{\phi}{\log \left( \frac{2}{\pi} \, (\phi-t) \right) \, dt} + 
\int_{\phi}^{\tau_j}{\log \left( \frac{2}{\pi} \, (t-\phi) \right)\, dt} \cr 
& \geq (\tau_j-\tau_{j-1})(c_9 - \log n) \cr \endsplit$$
with an absolute constant $c_9 > 0$, and (3.9) follows.
While, if 
$$\phi \not\in \left[\tau_{j-1} + \frac{c_8}{2n},\tau_j - \frac{c_8}{2n}\right]\,,$$ 
then Lemma 3.13 implies that there is an absolute constant  
$c_{10} > 0$ such that
$$\min_{t \in [\tau_{j-1},\tau_j]}|e^{it}-w| \geq c_{10}/n\,,$$
and hence 
$$\split \int_{\tau_{j-1}}^{\tau_j}{\log|e^{it}-w| \, dt} & \geq \int_{\tau_{j-1}}^{\tau_j}{\log(c_{10}/n) \, dt} \cr 
& \geq (\tau_j-\tau_{j-1})(-\log n - c_{11}) \cr \endsplit$$
with an absolute constant $c_{11} > 0$, and (3.9) follows again. 
Combining (3.8) and (3.9) and recalling that $g_w(t) :=  \log{|e^{it}-w|}$ and $\tau_j-\tau_{j-1} \leq 4\pi/n$, 
we obtain the inequality of the lemma.
\qed \enddemo

\proclaim{Lemma 3.14}
There is an absolute constant $c_{12} > 0$ such that 
$$\int_{\tau_{j-1}}^{\tau_j}{\log \left( \frac{R_n(t)}{n} \right)} \geq -\frac{c_{12}}{n}$$
for every $j \in \{2,3,\ldots,m+1\}$.
\endproclaim

\demo{Proof of Lemma 3.14}
Let, as before, $w_{\nu}, \nu=1,2,\ldots,n-1$, denote the zeros of $P_k$. Recall that $|w_\nu| \leq 2$ 
for each $\nu=1,2,\ldots,n-1$. We define the annuli
$$E_{j,q} := D(a_j,2^{q+3}\pi/n) \setminus D(a_j,2^{q+2}\pi/n)\,, \qquad q=1,2,\ldots,$$
and the disk
$$E_{j,0} := D(a_j,8\pi/n)\,.$$
Observe that the sets $E_{j,q}$ are pairwise disjoint and 
$${\Bbb C} = \bigcup_{j=0}^\infty{E_{j,q}}\,.$$
By Lemmas 3.2 and 3.3 there is an absolute constant $c_1>0$ (depending only on the explicitly given value of $\gamma$) such that 
$E_{j,q}$ contains at most $c_1n(2^{q+3}\pi/n+1)$ zeros of $P_k$ and $E_{j,0}$ contains at most $c_1(8\pi+1)$ zeros of $P_k$.
Hence Lemmas 3.11 and Lemma 3.12 give that
$$\split & \int_{\tau_{j-1}}^{\tau_j}{\log(R_n(t)) \, dt} - \frac 12 (\log(R_n(\tau_j)) - \log(R_n(\tau_{j-1})))(\tau_j - \tau_{j-1}) \cr
= & \sum_{\nu=1}^{n-1}{\left( \int_{\tau_{j-1}}^{\tau_j}{g_{w_\nu}(t) \, dt} - 
\frac 12 (g_{w_\nu}(\tau_j) - g_{w_\nu}(\tau_{j-1}))(\tau_j - \tau_{j-1}) \right)} \cr 
= & \sum_{q=0}^{\infty}{\sum_{w_{\nu} \in E_q}{ \left( \int_{\tau_{j-1}}^{\tau_j}{g_{w_\nu}(t) \, dt} - 
\frac 12 (g_{w_\nu}(\tau_j)-g_{w_\nu}(\tau_{j-1}))(\tau_j-\tau_{j-1}) \right)}} \cr
= & \sum_{q=0}^0{} + \sum_{q=1}^{\infty}{} 
\geq c_1(8\pi+1)\frac{-c_7}{n} + \sum_{q=1}^{\infty}{(c_1(n(2^{q+3}\pi/n)+1)) \frac{-c_6}{n^3(2^{q+2}/n)^2}} \cr
\geq & c_1(8\pi+1)\frac{-c_7}{n} - \sum_{q=1}^{\infty}{\frac{c_1c_6}{2^qn}} \cr 
\geq & -c_{12}/n \cr \endsplit$$ 
with an absolute constant $c_{12} > 0$. 
Now recall that $R_n(\tau_{j-1}) \geq 2\gamma n$ and $R_n(\tau_j) \geq 2\gamma n$, and the result follows. 
\qed \enddemo

Now we are ready to prove Lemma 3.8. 

\demo{Proof of Lemma 3.8}
Given $\varepsilon \in (0,1)$, let 
$$I_{n,\varepsilon} := \left\{ j \in \{2,3,\ldots,m+1\}: \min_{t \in [\tau_{j-1},\tau_j]}{R_n(t)} < \varepsilon \right\}\,,$$
and let
$$J_{n,\varepsilon} :=  \bigcup_{j \in I_{n,\varepsilon}}{[\tau_{j-1},\tau_j]} \,.$$
Using that $0 \leq R_n(t) \leq 2n$ for every $t \in K$, and then using Lemmas 3.4 and 3.14, we get  
$$\split \int_{B_{n,\varepsilon}}{\log \frac{R_n(t)}{2n} \, dt} & \geq \int_{J_{n,\varepsilon}}{\log \frac{R_n(t)}{2n} \, dt} 
= \sum_{j \in I_{n,\varepsilon}}{\int_{\tau_{j-1}}^{\tau_j}{\log \frac{R_n(t)}{2n} \, dt}} \cr 
& \geq c_3 n \varepsilon^{1/2}(-c_{12}/n) \geq -c_4 \varepsilon^{1/2} \cr \endsplit$$
for every sufficiently large $n \geq n_\varepsilon$, where $c_4 = c_3c_{12} > 0$, and the lemma is proved.
\qed \enddemo

\head 4. Proof of the Theorem \endhead

\demo{Proof of Theorem 2.1}
It follows from (1.2) immediately that 
$$\lim_{n \rightarrow \infty}{\frac{M_0(P_k)}{n^{1/2}}} = \lim_{n \rightarrow \infty}{\frac{M_0(Q_n)}{n^{1/2}}}\,,$$
so it is sufficient to prove the asymptotic formula only for $M_0(P_k)$.

Let $\varepsilon \in (0,1)$ be fixed. By Lemma 3.7 we have
$$\lim_{n \rightarrow \infty}{\frac{1}{2\pi} \int_{A_{n,\varepsilon}}{\log \frac{R_n(t)}{2n}\, dt}} = \int_\varepsilon^1{\log x\, dx}\,.$$
while it follows from Lemma 3.8 and the inequalities $0 \leq R_n(t) \leq 2n$ that there is an absolute constant $c_4 > 0$ such that 
$$-c_4\varepsilon^{1/2} \leq \frac{1}{2\pi} \int_{B_{n,\varepsilon}}{\log \frac{R_n(t)}{2n}\, dt} \leq 0$$
for every sufficiently large $n \geq n_\varepsilon$. As $K$ is the disjoint union of $A_{n,\varepsilon}$ and $B_{n,\varepsilon}$, 
we have
$$\limsup_{n \rightarrow \infty}{\frac{1}{2\pi} \int_K{\log \frac{R_n(t)}{2n}\, dt}} \leq \int_\varepsilon^1{\log x\, dx} \tag 4.1$$
and 
$$\liminf_{n \rightarrow \infty}{\frac{1}{2\pi} \int_K{\log \frac{R_n(t)}{2n}\, dt}} \geq \int_\varepsilon^1{\log x\, dx} 
- c_4\varepsilon^{1/2}\,. \tag 4.2$$
As (4.1) and (4.2) hold for an arbitrary $\varepsilon \in (0,1)$, it follows that
$$\int_0^1{\log x\, dx} \leq \liminf_{n \rightarrow \infty}{\frac{1}{2\pi} \int_K{\log \frac{R_n(t)}{2n}\, dt}} 
\leq \limsup_{n \rightarrow \infty}{\frac{1}{2\pi} \int_K{\log \frac{R_n(t)}{2n}\, dt}} \leq \int_0^1{\log x\, dx}\,,$$
and hence
$$\lim_{n \rightarrow \infty}{\frac{1}{2\pi} \int_K{\log \frac{R_n(t)}{2n}\, dt}} =-1\,.$$
Hence, recalling that 
$$R_n(t) := |P_k(e^{it}|^2, \qquad t \in K\,,$$
we obtain
$$\split \lim_{n \rightarrow \infty}{\frac{1}{2\pi} \int_K{\log \frac{|P_k(e^{it})|}{(2n)^{1/2}}\, dt}}
= & \lim_{n \rightarrow \infty}{\frac{1}{2\pi} \int_K{\log \left( \frac{R_n(t)}{2n} \right)^{1/2} \, dt}} \cr 
= & \lim_{n \rightarrow \infty}{\frac{1}{2\pi} \frac 12 \, \int_K{\log \frac{R_n(t)}{2n}  \, dt}} = -1/2\,. \cr \endsplit$$
Hence
$$\split \lim_{n \rightarrow \infty}{\frac{M_0(P_k)}{(2n)^{1/2}}} 
= & \lim_{n \rightarrow \infty}{\exp \left( \frac{1}{2\pi} \int_K{\log \frac{|P_k(e^{it})|}{(2n)^{1/2}}\, dt} \right)} \cr  
= & \exp \left( \lim_{n \rightarrow \infty}{\frac{1}{2\pi} \int_K{\log \frac{|P_k(e^{it})|}{(2n)^{1/2}}\, dt}} \right) = \exp(-1/2) \cr \endsplit$$
which is the asymptotic formula for $M_0(P_k)$ stated in the theorem.
\qed \enddemo

\head 5. The Mahler measure of the Fekete polynomials \endhead

For a prime $p$ the $p$-th Fekete polynomial is defined as
$$f_p(z) := \sum_{k=1}^{p-1}{\left( \frac kp \right)z^k}\,,$$
where
$$\left( \frac kp \right) =
\cases
1, \quad \text{if \enskip} x^2 \equiv k \enskip (\text {mod\,}p) \enskip
\text{for an} \enskip x \not\equiv 0 \enskip (\text {mod\,}p)\,,
\\
0, \quad \text{if \enskip} p \enskip \text{divides} \enskip k\,,
\\
-1, \quad \text{otherwise}
\endcases$$
is the usual Legendre symbol. Since $f_p$ has constant coefficient $0$, it is not
a Littlewood polynomial, but $g_p$ defined by $g_p(z) := f_p(z)/z$ is a Littlewood
polynomial of degree $p-2$, and has the same Mahler measure as $f_p$. Fekete polynomials 
are examined in detail in [B-02], [CG-00], [E-11], [E-12], [EL-07], and [M-80]. 
In [CE-15a] and [CE-15b] the authors examined the maximal size of the Mahler measure of 
sums of $n$ monomials on the unit circle as well as on subarcs of the unit circles. In 
the constructions appearing in [CE-15a] properties of the Fekete polynomials $f_p$ 
turned out to be quite useful. Montgomery [M-80] proved the following fundamental result.

\proclaim{Theorem 5.1}
There are absolute constants $c_{13} > 0$ and $c_{14} > 0$ such that
$$c_{13} \sqrt{p} \log \log p \leq \max_{z \in \partial D}{|f_p(z)|} \leq c_{14} \sqrt{p} \log p\,.$$
\endproclaim

In [E-07] we proved the following result.

\proclaim{Theorem 5.2} For every $\varepsilon > 0$ there is
a constant $c_{\varepsilon}$ such that
$$M_0(f_p) \geq \left(\frac 12 - \varepsilon \right)\sqrt{p}$$
for all primes $p \geq c_{\varepsilon}$.
\endproclaim

In [E-17] the factor $\left(\frac 12 - \varepsilon \right)$ in Theorem 1.2 has been 
improved to to an absolute constant $c > 1/2$. Namely we prove the following.

\proclaim{Theorem 5.3} There is an absolute constant $c > 1/2$ such that
$$M_0(f_p) \geq c\sqrt{p}$$
for all sufficiently large primes.
\endproclaim

The determine the asymptotic size of the Mahler measure $M_0(f_p)$ of the Fekete polynomials $f_p$ seems 
to be beyond reach at the moment. Not even a (published or unpublished) conjecture seems to be known.

\head 6. Acknowledgement \endhead
The author thanks Stephen Choi and Bahman Saffari for checking the details of the proof 
in this paper and for their suggestions to make the paper more readable.

\enddocument